\documentclass[12pt]{article}
\usepackage{amsmath, amssymb, amscd}
\usepackage[matrix,arrow,curve]{xy}

\newcommand{\Var}{{\cal{V}_{\mathbb{C}}}}

\newcommand{\St}{{\rm Stck}_{\mathbb{C}}}

\def\1{\underline{1}}

\def\AA{{\mathbb A}}

\def\bL{{\mathbb L}}

\def\C{{\mathbb C}}

\def\CP{\mathbb C\mathbb P}

\def\Gr{\mbox{\bf Gr}}

\newtheorem{theorem}{Theorem}

\newtheorem{lemma}{Lemma}

\newenvironment{definition}
{\smallskip\noindent{\bf Definition\/}:}{\smallskip\par}

\newenvironment{proof}
{\noindent{\bf Proof\/}.}{{ $\square$}\smallskip\par}

\title{On piecewise isomorphism of some varieties.
}
\author{S.M.~Gusein-Zade \thanks{Partially supported by the grants
RFBR-10-01-00678,
NSh-8462.2010.1.
Address: Moscow State University, Faculty
of Mathematics and Mechanics, Moscow, GSP-1, 119991, Russia. E-mail:
sabir\symbol{'100}mccme.ru} \and I.~Luengo \thanks{The last two authors are partially
supported by the grants MTM2010-21740-C02-01 and Grupo Singular CCG07-UCM/ESP-2695-921020. Address: University
Complutense de Madrid, Dept. of Algebra, Madrid, 28040, Spain.
E-mail: iluengo\symbol{'100}mat.ucm.es, amelle\symbol{'100}mat.ucm.es} \and
A.~Melle--Hern\'andez \thanks{Address: Instituto de Ciencias Matem\'aticas 
Complutense-CSIC-Aut\'onoma-Carlos III, Spain.
}} 
 \date{}
\begin{document}

\maketitle

\begin{abstract}
Two quasi-projective varieties  are called {\em piecewise
isomorphic} if they can be stratified into  pairwise isomorphic strata. We show that the $m$-th symmetric power $S^m(\C^n)$
of the complex affine space $\C^n$ is piecewise isomorphic to $\C^{mn}$ and the $m$-th symmetric power $S^m(\CP^\infty)$
of the infinite dimensional complex projective space is  piecewise isomorphic to the infinite dimensional Grassmannian $\Gr(m,\infty)$.

\end{abstract}


Let $K_0(\Var)$ be the Grothendieck ring of complex
quasi-projective varieties. This is the Abelian group generated by
the classes $[X]$ of all complex quasi-projective varieties $X$
modulo the relations:
\begin{enumerate}
\item[1)] $[X]=[Y]$ for isomorphic $X$ and $Y$;
\item[2)] $[X]=[Y]+[X\setminus Y]$ when $Y$ is a Zariski closed subvariety of $X$.
\end{enumerate}
The multiplication in $K_0(\Var)$ is defined by the Cartesian
product of varieties: $[X_1]\cdot [X_2]=[X_1\times X_2]$. The
class $[\AA^1_\C]\in K_0(\Var)$ of the complex affine line is
denoted by $\bL$.

\begin{definition}
Quasi-projective varieties $X$ and $Y$ are called {\em piecewise
isomorphic} if there exist decompositions
$X=\coprod\limits_{i=1}^{s}X_i$ and
$Y=\coprod\limits_{i=1}^{s}Y_i$ of $X$ and $Y$ into (Zariski)
locally closed subsets such that $X_i$ and $Y_i$ are isomorphic
for $i=1, \ldots, s$.
\end{definition}

If the varieties $X$ and $Y$ are piecewise isomorphic, their
classes $[X]$ and $[Y]$ in the Grothendieck ring $K_0(\Var)$
coincide. There exists the conjecture (or at least the
corresponding question) that the opposite also holds: if
$[X]=[Y]$, then $X$ and $Y$ are piecewise isomorphic: see
\cite{LL, LS}.

It is well-known that the $m$-th symmetric power $S^m\C^n$ of the
affine space $\C^n$ is birationally equivalent to $\C^{mn}$: see
e.g. \cite{GKZ}. An explicit birational isomorphism between
$S^m\C^n$ and $\C^{mn}$ was constructed in \cite{BR}. Moreover the
class $[S^m\C^n]$ of the variety $S^m\C^n$ in the Grothendieck
ring $K_0(\Var)$ of complex quasi-projective varieties is equal to
the class $[\C^{mn}]={\bL}^{mn}$: see e.g. \cite{G, GLM04}. The conjecture formulated
above means that the varieties $S^m\C^n$ and $\C^{mn}$ are
piecewise isomorphic. This is well-known for $n=1$. Moreover $S^m\C$ and $\C^{m}$
are isomorphic. 
The fact that indeed $S^m\C^n$ and $\C^{mn}$  are  piecewise
isomorphic seems to (or must) be known to specialists.
Moreover proofs are essentially contained in \cite{G} (Lemma 4.4 proved by Burt Totaro) and \cite{GLM04} (Statement 3).
However this fact is not explicitly reflected in the literature. Here we give a proof
of this statement.

In \cite{stacks}, it was shown that the Kapranov zeta function
$\zeta_{B\C^*}(T)$ of the classifying stack $B\C^*=BGL(1)$ is
equal to
$$
1+\sum_{m=1}^\infty
\frac{\bL^{m^2-m}}{(\bL^m-\bL^{m-1})(\bL^m-\bL^{m-2})
\ldots(\bL^m-1)} T^m.
$$
Unrigorously speaking this can be interpreted as the class $[S^m
B\C^*]$ of the ``$m$-th symmetric power'' of the classifying stack
$B\C^*$ in the Grothendieck ring $K_0(\St)$ of algebraic stacks of
finite type over $\C$ is equal to $\bL^{m^2-m}$ times the class
$[BGL(m)]=1/(\bL^m-\bL^{m-1})(\bL^m-\bL^{m-2})\ldots(\bL^m-1)$ of
the classifying stack $BGL(m)$. The natural topological analogues
of the classifying stacks $B\C^*$ and $BGL(m)$ are the
infinite-dimensional projective space $\CP^{\infty}$ and the
infinite Grassmannian 
$\Gr(m,\infty)$. We
show that the $m$-th symmetric power $S^m \CP^{\infty}$ of
$\CP^{\infty}$ and $\Gr(m,\infty)$ are piecewise isomorphic in a
natural sense.

\begin{theorem}\label{prop1}
The varieties $S^m\C^n$ and $\C^{mn}$ are piecewise isomorphic.
\end{theorem}

\begin{proof}
The proofs which we know in any case are not explicit, we do not know the neccesary partitions of 
$S^m\C^n$ and  $\C^{mn}$. Therefore we prefer to use the language of power structure
over the Grothendieck semiring $S_0({\rm{Var}}_{\C} )$ of complex quasi-projective varieties
invented in \cite{GLM04}. This language sometimes permits to substitute
somewhat envolved combinatorial considerations by short computations
(or even to avoid them at all, as it was made in \cite{GLM06}).
Since the majority of statements in \cite{GLM04} (including those which 
could be used to prove Theorem~1) are formulated and proved in
the Grothendieck ring $K_0(\Var)$  of complex quasi-projective varieties,
we repeat a part of the construction in the appropriate setting. 

The Grothendieck semiring $S_0({\rm{Var}}_{\C})$ of complex quasi-projective
varieties is the semigroup generated by isomorphism classes $\{X\}$
of such varieties modulo the relation $\{X\}=\{X-Y\}+\{Y\}$ for a
Zariski closed subvariety $Y\subset X$. The
multiplication is defined by the Cartesian product of varieties:
$\{X_1\}\cdot \{X_2\}=\{X_1\times X_2\}.$
Classes $\{X\}$ and $\{Y\}$ of two varieties $X$ and $Y$ in $S_0({\rm{Var}}_{\C} )$
are equal if and only if $X$ and $Y$ are piecewise isomorphic.
Let $\bL\in S_0({\rm{Var}}_{\C})$
be the class of the affine line. If $\pi:E \to B$ is a Zariski locally trivial fibre bundle 
with fibre $F$, one has $\{E\}=\{F\}\cdot \{B\}$. For example if $\pi:E \to B$ is a Zariski locally trivial 
vector bundle of rank $s$, one has $\{E\}=\bL^s \{B\}$.

A {\em power structure} over a semiring $R$ is a map
$\left(1+T\cdot R[[T]]\right)\times {R} \to 1+T\cdot R[[T]]$:
$(A(T),m)\mapsto \left(A(T)\right)^{m}$,
which possesses the properties:
\begin{enumerate}
\item $\left(A(T)\right)^0=1$,
\item $\left(A(T)\right)^1=A(T)$,
\item $\left(A(T)\cdot B(T)\right)^{m}=\left(A(T)\right)^{m}\cdot
\left(B(T)\right)^{m}$,
\item $\left(A(T)\right)^{m+n}=\left(A(T)\right)^{m}\cdot
\left(A(T)\right)^{n}$,
\item $\left(A(T)\right)^{mn}=\left(\left(A(T)\right)^{n}\right)^{m}$,
\item $(1+T)^m=1+mT+$ terms of higher degree,
\item $\left(A(T^\ell)\right)^m =
\left(A(T)\right)^m\raisebox{-0.5ex}{$\vert$}{}_{T\mapsto T^\ell}$, $\ell\geq 1$.
\end{enumerate}

In \cite{GLM04}, there was defined a power structure over the Grothendieck semiring $S_0({\rm{Var}}_{\C} )$.
Namely, for $A(T)=1+\{A_1\}\,T+\{A_2\}\,T^2+\ldots$ and $\{M\} \in S_0({\rm{Var}}_{\C})$, the series $\left(A(T)\right)^{\{M\}}$ 
is defined as 
\begin{equation} \label{eq00}
1+\sum_{k=1}^\infty
\left(\sum_{\{k_i\}:\sum ik_i=k}
\left\{
\left( (
(\prod_i M^{k_i}
)
\setminus\Delta
)
\times\prod_i A_i^{k_i}\right)/\prod_i S_{k_i}
\right\}
\right)
\cdot T^k,
\end{equation}
where $\Delta$ is the "large diagonal" in $M^{\Sigma k_i}=\prod\limits_i M^{k_i}$
which consists of
$(\sum k_i)$-tuples of points of $M$ with at least two coinciding ones,
the group $S_{k_i}$ of permutations on $k_i$ elements acts by permuting corresponding $k_i$ factors in
$\prod\limits_i M^{k_i}\supset (\prod\limits_i M^{k_i})\setminus\Delta$
and the spaces $A_i$ simultaneously. The action of the group
$\prod\limits_i S_{k_i}$ on $(\prod\limits_i M^{k_i})\setminus\Delta$ is free.
The properties 1--7 are proved in \cite[Theorem 1]{GLM04}.

Special role is played by the Kapranov zeta function in the Grothendieck
semiring $S_0({\rm{Var}}_{\C} )$:
$$
\zeta_{\{M\}}(T):=1+\sum_{k=1}^\infty \{S^kM\} T^k\, ,
$$
where $S^kM$ is the $k$-th symmetric power $M^k/S_k$ of the variety $M$.
In terms of the power structure one has $\zeta_{\{M\}}(T)=(1+T+T^2+\ldots)^{\{M\}}$.
Theorem~1 is equivalent to the fact that 
\begin{equation}\label{ouraim}
\zeta_{\bL^m}(T)=(1+\sum_{i=1}^\infty \bL^{im} T^i). 
\end{equation}

\begin{lemma}\label{lemma}
Let $A_i$ and $M$ be complex quasi-projective varieties,
$A(T)=1+\{A_1\}T+\{A_2\}T^2+\ldots$. Then, for any integer $s\ge0$,
\begin{equation} \label{from-lemma}
\left(A(\bL^s T)\right)^{\{M\}} = \left(A(T)^{\{M\}}\right)
\mbox{\raisebox{-0.5ex}{$\vert$}}{}_{T\mapsto {\bL^s T}}.
\end{equation}
\end{lemma}

\begin{proof}
The coefficient at the monomial $T^k$ in the power series
$\left(A(T)\right)^{\{M\}}$
is a sum of the classes of varieties of the form
$$
V=\left( (
(\prod_i M^{k_i}
)
\setminus\Delta
)
\times\prod_i A_i^{k_i}\right)/\prod_i S_{k_i}
$$
with $\sum ik_i=k.$ The corresponding summand $\{\widetilde V\}$ in the coefficient
at the monomial $T^k$ in the power series $\left(A(\bL^s T)\right)^{\{M\}}$
has the form
$$\widetilde V=\left((
(\prod_i M^{k_i}) \setminus \Delta
)
\times\prod_i (\bL^{si}\times A_i)^{k_i}\right)/\prod_iS_{k_i}.
$$
The natural map $\widetilde V \to V$
 is a Zariski locally
trivial vector bundle of rank $sk$: see e.g.\cite[Section 7, Proposition 7]{Mumf}. 
This implies that $\{\widetilde V\}=\bL^{sk}\cdot\{V\}.$
\end{proof}

One has $\zeta_{\bL}(T)=(1+\bL T+ \bL^2 T^2+\ldots)$. For all $A_i$ being points, 
i.e. $\{A_i\}=1$, one gets
\begin{eqnarray*}
\zeta_{\bL \{M\}}(T)&=&(1+T+T^2+\ldots)^{\bL \{M\}}=({(1+T+T^2+\ldots)^{\bL}})^{\{M\}}\\
&=&
(1+\bL T+ \bL^2 T^2+\ldots)^{\{M\}}.
\end{eqnarray*}
Equation~(\ref{from-lemma}) implies that 
$$
\zeta_{\bL \{M\}}(T)=(1+\bL T+ \bL^2 T^2+\ldots)^{\{M\}}=\zeta_{\{M\}}(\bL T).
$$
Assuming (\ref{ouraim}) holds for $m<m_0$ and applying the equation above
to $m=m_0-1$ one gets 
\begin{eqnarray*}
\zeta_{\bL^{m_0}}(T)&=& \zeta_{\bL^{m_0-1}}(\bL T)=(1+\bL^{m_0-1} T+ \bL^{2(m_0-1)} T^2+\ldots)\mbox{\raisebox{-0.5ex}{$\vert$}}{}_{T\mapsto {\bL T}}\\
&=& (1+\bL^{m_0} T+ \bL^{2m_0} T^2+\ldots).
\end{eqnarray*}
This gives the proof.
\end{proof}

Let $\CP^{\infty}=\varinjlim\CP^N$ be the infinite dimensional
projective space and let $\Gr(m,\infty)=\varinjlim \Gr(m,N)$ be the
infinite dimensional Grassmannian. (In the both cases the
inductive limit is with respect to the natural sequence of
inclusions. The spaces $\CP^{\infty}$ and $\Gr(m,\infty)$ are,
in the topological sense, classifying spaces for the groups
$\C^*=GL(1;\C)$ and $GL(m;\C)$ respectively.) The symmetric power
$S^m \CP^{\infty}$ is the inductive limit of the quasi-projective
varieties $S^m \CP^{N}$. For a sequence $X_1\subset X_2\subset X_3
\subset \ldots$ of quasi-projective varieties, let $X=\varinjlim
X_i(=\bigcup\limits_i X_i)$ be its (inductive) limit. A \emph{partition
of the space $X$ compatible with the filtration} $\{X_i\}$ is a
representation of $X$ as a disjoint union $\coprod\limits_j Z_j$
of (not more than) countably many quasi-projective varieties $Z_j$
such that each $X_i$ is the union of a subset of the strata $Z_j$
and each $Z_j$ is a Zariski locally closed subset in the
corresponding $X_i$.

\begin{theorem}\label{prop2}
The spaces $S^m \CP^{\infty}$ and $\Gr(m,\infty)$ are piecewise
isomorphic in the sense that there exist partitions $S^m
\CP^{\infty}=\coprod\limits_j U_j$ and
$\Gr(m,\infty)=\coprod\limits_j V_j$ into pairwise isomorphic
quasi-projective varieties $U_j$ and $V_j$ ($U_j\cong V_j$)
compatible with the filtrations $\{S^m \CP^{N}\}_N$ and
$\{\Gr(m,N)\}_N$.
\end{theorem}

\begin{proof}
The natural partition of $\Gr(m,N)$ consists of the Schubert cells
corresponding to the flag $\{0\}\subset \C^1\subset \C^2\subset\ldots$: see e.g \cite[\S 5.4]{FF}.
Each Schubert cell is a locally closed subvariety of $\Gr(m,N)$
isomorphic to the complex affine space of certain dimension. This partition
is compatible with the inclusion $\Gr(m,N)\subset\Gr(m,N+1)$ and
therefore gives a partition of $\Gr(m,\infty)$. The number of cells of
dimension $n$ in $\Gr(m,\infty)$ is equal to the number of partitions of
$n$ into summands not exceeding $m$.

Since $\CP^{\infty}=\C^0\coprod \C^1\coprod \C^2\coprod\ldots$ and
$S^p(A\coprod B)=\coprod\limits_{i=0}^p S^iA\times S^{p-i}B$,
one has 
$$
S^m \CP^{\infty}=\coprod_{\{i_0,i_1,i_2,\ldots\}:i_0+i_1+i_2+\ldots=m}
\prod_{j} S^{i_j}\C^j=\coprod_{\{i_1,i_2,\ldots\}:i_1+i_2+\ldots\le m}
\prod_j S^{i_j}\C^j,
$$
where $i_j$ are non-negative integers.
This partition is compatible with the natural filtration
$\{\CP^0\}\subset \CP^1\subset \CP^2\subset\ldots$.
The number of parts of dimension $n$ is equal to the
number of sequences $\{i_1,i_2,\ldots\}$ such that $\sum\limits_j i_j\le m$,
$\sum\limits_j i_j j=n$. Thus it coincides with the number of partitions of $n$
into not more than $m$ summands and is equal to the number of
$n$-dimensional Schubert cells in the partition of $\Gr(m,\infty)$.
Due to Proposition~1 each part $\prod\limits_j S^{i_j}\C^j$ is piecewise
isomorphic to the complex affine space of the same dimension.
This concludes the proof.
\end{proof}
It would be interesting to find explicit piecwise isomorphisms between
the spaces in Theorems~1 and~2.


\begin{thebibliography}{15}

\bibitem{BR} V.M.~Buchstaber, E.G.~Rees. \emph{Manifolds of polysymmetric
polynomials. Classical problems, contemporary applications}. In:
Proceedings of the Conference devoted to the 10th anniversary of
RFBR, Fizmatlit, Moscow, 2004, 129--145 (in Russian).

\bibitem{FF} A.T.~Fomenko, D.B.~Fuks.
\emph{A course in homotopic topology}. Nauka, Moscow, 1989. 496 pp.  (in Russian).

\bibitem{GKZ} I.M.~Gelfand, M.M.~Kapranov, A.V.~Zelevinsky. \emph{Discriminants,
resultants, and multidimensional determinants}. Birkhauser, Boston,
1994.

\bibitem{G} L.~G\"ottsche. \emph{On the motive of the Hilbert scheme
of points on a surface}. Math. Res. Lett., v.8, no.5--6, 613--627 (2001).

\bibitem{GLM04} S.M.~Gusein-Zade, I.~Luengo, A.~Melle-Hern\'andez.
\emph{A power structure over the Grothendieck ring of varieties}.
Math. Res. Lett., v.11, no.1, 49--57 (2004).

\bibitem{GLM06} S.M.~Gusein-Zade, I.~Luengo, A.~Melle-Hern\'andez.
\emph{Power structure over the Grothendieck ring of varieties and
generating series of Hilbert schemes of points}. Michigan Math. J.,
v.54, 353--359 (2006).

\bibitem{stacks} S.M.~Gusein-Zade, I.~Luengo, A.~Melle-Hern\'andez.
\emph{On the pre-lambda ring structure on the Grothendieck ring of
stacks and on the power structures over it}. ArXiv: 1008.5063.

\bibitem{LL} M.~Larsen, V.~Lunts. \emph{Motivic measures and stable
birational geometry}. Moscow Math. J., v.3, no.1, 85--95 (2003).

\bibitem{LS} Q.~Liu, 
J.~Sebag. \emph{The Grothendieck ring of varieties and piecewise
isomorphisms}. Mathematische Zeitschrift, v.265, no.2, 321--342 (2010).

\bibitem{Mumf} D.~Mumford. \emph{Abelian varieties}. Tata Institute of Fundamental Research Studies in 
Mathematics, Bombay, No. 5;
Oxford University Press, London, 1970, viii+242 pp. 

\end{thebibliography}
\end{document}